\documentclass[12pt]{amsart}
\usepackage{amsmath,amssymb}

\newtheorem{thm}{Theorem}
\newtheorem{lem}[thm]{Lemma}
\newtheorem*{prop2}{Proposition}
\newtheorem{cor}[thm]{Corollary}

\newcommand{\Z}{{\mathbf Z}}
\newcommand{\Q}{{\mathbf Q}}

\newcommand{\R}{{\mathbf R}}

\newcommand{\ord}{\operatorname{ord}}
\newcommand{\trace}{\operatorname{tr}}
\newcommand{\gal}{\operatorname{Gal}}
\newcommand{\disc}{\operatorname{disc}}

\newcommand{\A}{A}

\newcommand{\unit}{\epsilon}
\newcommand{\unitp}{\overline{\epsilon}}
\newcommand{\OO}{\mathfrak{O}}

\newcommand{\OK}{\mathfrak{O}_K}

\title{On the order of unimodular matrices modulo integers}

\author{P\"ar Kurlberg}
\address{
Department of Mathematics\\ 
Chalmers University of Techology\\
SE-412 96 Gothenburg  \\
Sweden}
\urladdr{www.math.chalmers.se/\~{ }kurlberg}
\email{kurlberg@math.chalmers.se}
\thanks{Author supported in part by the National Science Foundation
  (DMS 0071503).}

\begin{document}

\begin{center}
  PRELIMINARY VERSION
\end{center}

\begin{abstract}
Assuming the Generalized Riemann Hypothesis, we prove the following:  
If  $b$ is an integer greater than one, then the
multiplicative order of $b$ 
modulo $N$ is larger than $N^{1-\epsilon}$ for all $N$ in a density
one subset of the integers.
If $A$ is a hyperbolic unimodular matrix with integer coefficients,
then the order of $A$ modulo $p$ is greater than $p^{1-\epsilon}$ for
all $p$ in a density one subset of the primes.  Moreover,  the 
order of $A$ modulo $N$ is greater than $N^{1-\epsilon}$ for all $N$ 
in a density one subset of the integers.
\end{abstract}

\maketitle




\section{Introduction}
Given an integer $b$ and a prime $p$ such that $p {\not|} b$, let
$\ord_p(b)$ be the multiplicative 
order of $b$ modulo $p$.  In other words, $\ord_p(b)$ is the smallest
non negative integer $k$ such that $b^k \equiv 1 \mod p$.  Clearly
$\ord_p(b) \leq p-1$, and if the order is maximal, $b$ is said to be a
primitive root modulo $p$. 
%
Artin conjectured (see the preface in \cite{artin-collected}) that if 
$b \in \Z$ is not 
a square, then $b$ is a primitive root for a positive 
proportion\footnote{The constant is given by an Euler product that
  depends on $b$.} of 
the primes.

What about the ``typical'' behaviour of $\ord_p(b)$?  For
instance, are there good lower bounds on $\ord_p(b)$ that hold for a
full 
density 
subset of the primes?  In \cite{erdos-murty}, Erd{\H{o}}s
and Murty proved that if $b \neq 0,\pm 1$, then there exists a
$\delta>0$ so that $\ord_p(b)$ is at least $p^{1/2} \exp( (\log
p)^\delta)$ for a full density subset of the primes.
However, we expect the typical order to be much larger.  In
\cite{hooley-artin} Hooley proved that the Generalized Riemann
Hypothesis (GRH) implies Artin's conjecture.
Moreover, if $f : \R^+ \to \R^+$ is an increasing function tending
to infinity, Erd{\H{o}}s and Murty
showed
\cite{erdos-murty} that GRH implies that the order of $b$ modulo $p$
is greater 
than $p/f(p)$ for full density subset of the primes.

It is also interesting to consider lower bounds for $\ord_N(b)$ where
$N$ is an integer.  It is easy to see that $\ord_N(b)$ can be as small
as $\log N$ infinetely often (take $N = b^k-1$), but we expect that the
typical order to be quite large.  Assuming GRH, we can prove that
the lower bound $\ord_N(b) \gg N^{1-\epsilon}$ holds for most
integers.
\begin{thm}
\label{t:ord-b-all-N}
Let $b \neq 0,\pm 1$ be an integer.  Assuming GRH, the number of $N\leq
x$ such that $\ord_N(b) \ll N^{1-\epsilon}$ is $o(x)$. That is, the
set of integers $N$ such that $\ord_N(b) \gg N^{1-\epsilon}$ has
density one.
\end{thm}

However, the main focus of this paper is to investigate a related
question, 
namely lower bounds on the order of unimodular matrices modulo $N \in
\Z$.  That is, if $A \in SL_2(\Z)$, what can be said about lower
bounds for $\ord_N(A)$, the order of $A$ modulo $N$, that hold for
most $N$?
It is a  natural generalization of the previous questions, but
our main motivation comes from mathematical physics (quantum chaos):
In \cite{cat2} Rudnick and I proved that if $A$ is
hyperbolic\footnote{$A$ is hyperbolic if $|\trace(A)|>2$.}, then
quantum ergodicity for toral automorphisms follows from $\ord_N(A)$
being slightly larger than $N^{1/2}$, and we then showed that
this condition does hold for a full density subset of the integers.
%

Again, we expect that the typical order is much larger.
In order to give lower bounds on $\ord_N(A)$, it is
essential to have 
good lower bounds on $\ord_p(A)$ for $p$ prime:
%
\begin{thm}
\label{t:hooley++}
Let $A \in SL_2(\Z)$ be hyperbolic, and let $f : \R^+ \to \R^+$ be an
increasing function tending to 
infinity slower than $\log x$.  Assuming GRH, there are at most
$O(\frac{x }{\log x 
  f(x)^{1-\epsilon} 
})$ primes $p \leq x$ such that $\ord_p(A) < p/f(p)$.  In particular, the
set of primes $p$ such that $\ord_p(\A) \geq 
p/f(p)$ has density one.
\end{thm}

Using this we  obtain an improved lower bound on $\ord_N(A)$
that is valid for most integers.
\begin{thm}
\label{t:all-N}
Let $A \in SL_2(\Z)$ be hyperbolic.   Assuming GRH, the number of
  $N\leq x$ such that $\ord_N(A) \ll 
  N^{1-\epsilon}$ is $o(x)$. That is, the set of integers $N$ such
  that  $\ord_N(A)  \gg   N^{1-\epsilon}$ has density one.
\end{thm}

{\em Remarks:}
If $A$ is elliptic ($|\trace(A)|<2$) then $A$ has finite order (in
fact, at most $6$).
If $A$ is parabolic ($|\trace(A)|=2$), then $\ord_p(A)=p$ unless
$A$ is congruent to the identity matrix modulo $p$, and 
hence there exists a constant $c_A>0$ so that $\ord_N(A) > c_A N$.
Apart from the application in mind, it is thus natural to  only
treat the hyperbolic case.

As far as unconditional results for primes go, we note that the
proof in \cite{erdos-murty} relies entirely on analyzing the divisor
structure of $p-1$, and we expect that their method should give a
similar lower bound on the order of $A$ modulo $p$.
An unconditional
lower bound of the form
\begin{equation}
\label{eq:order-lower-exp-bound}
\ord_p(b) \gg p^{\eta}
\end{equation}
for a full proportion of the primes and $\eta>1/2$ would be quite
interesting.  In this direction, Goldfeld proved
\cite{goldfeld-big-prime-divisor} that  if $\eta<3/5$, then
\eqref{eq:order-lower-exp-bound} holds for a positive, but not full,
proportion of the primes.


%

Clearly $\ord_p(A)$ is related to $\ord_p(\unit)$, where $\unit$ is
one of the eigenvalues of $A$.  Since $A$ is assumed to be hyperbolic,
$\unit$ is a power of a fundamental unit in a real quadratic field.
The question of densities of primes $p$ such that $\ord_p(\lambda)$ is
maximal, for $\lambda$ a fundamental unit in a real quadratic field,
does not seem to have received much attention until  quite
recently; in \cite{roskam-unit-artin} Roskam proved  that GRH implies that
the set of primes $p$ for which $\ord_p(\lambda)$ is maximal has
positive density.  (The work of Weinberger
\cite{cooke-weinberger-division-chains}, Cooke and Weinberger
\cite{weinberger-euclidean-rings} and Lenstra \cite{lenstra-artin}
does treat the case $\ord_p(\lambda) =p-1$, but not the case
$\ord_p(\lambda) =p+1$.)

\section{Preliminaries}
\subsection{Notation}

If $\OO_F$ is the ring of integers in a number field $F$, we let
$\zeta_F(s) = \sum_{\mathfrak{a} \subset \OO_F} N(\mathfrak{a})^{-s}$
denote the zeta function of $F$.  By GRH we mean that all nontrivial
zeroes of $\zeta_F(s)$ lie on the line $Re(s)=1/2$ for all number
fields $F$.  

Let $\unit$ be an eigenvalue of $\A$, satisfying the equation
\begin{equation}
\label{e:unit}
\unit^2-\trace(\A)\unit+\det(\A).
\end{equation}
Since $A$ is hyperbolic, $K= \Q(\unit)$ is a real quadratic field.
Let $\OK$ be the integers in $K$, and let $D_K$ be the discriminant of
$K$.  Since $A$ has determinant one, $\unit$ is a unit in
$\OK$. 
For $n \in \Z^{+}$ we let $\zeta_n = e^{2\pi i/n}$ be a primitive
$n$-th root of unity, and $\alpha_n = \unit^{1/n}$ be an $n$-th root
of $\epsilon$.  Further, with $Z_n = K(\zeta_n)$, $K_n=K(\zeta_n,
\alpha_n)$, and $L_n=K(\alpha_n)$, we let $\sigma_p$ denote the
Frobenius element in $\gal(K_n/\Q)$ associated with $p$. We let
$F_{p^k}$ denote the 
finite field with $p^k$ elements, and we let $F_{p^2}^1 \subset
F_{p^2}^\times$ be the norm one elements in $F_{p^2}$, i.e., the
kernel of the norm map from $F_{p^2}^\times$ to $F_{p}^\times$.
Let $\left<A\right>_p$ be the group generated by $A$ in
$SL_2(F_p)$. $\left<A\right>_p$ is 
contained in a maximal torus (of order $p-1$ or $p+1$), and we let
$i_p$ be the index of $\left<A\right>_p$ in this torus.
Finally, let $\pi(x) = |\{ p \leq x : \text{$p$ is prime}\}|$ be the
number of primes up to $x$.

\subsection{Kummer extensions and Frobenius elements}

We want to characterize primes $p$ such that $n|i_p$, and we can
relate this to primes splitting in certain Galois extensions as
follows:

Reduce equation \eqref{e:unit} modulo $p$ and
let $\unitp$ denote a solution to equation \eqref{e:unit} in $F_p$
or $F_{p^2}$. (Note that if $p$ does not ramify in $K$ then the order
of $\A$ modulo $p$ equals the order of $\unit$ modulo $p$.)
%
%
%
If $p$ splits in $K$ then $\unitp \in F_p$, and if $p$ is inert, then
$\unitp \in F_{p^2} \setminus F_p$. In the latter case, $\unitp \in
F_{p^2}^1$ since the norm one property is preserved when reducing
modulo $p$. Now, $F_p^\times$ and $F_{p^2}^1$ are cyclic groups of
order $p-1$ and $p+1$ respectively.  Thus, if $p$ splits in $K$ then
$\ord_p(\unit)|p-1$, whereas if $p$ is inert in $K$ then
$\ord_p(\unit)|p+1$.

\begin{lem}
\label{l:size-of-frob-conj-class}
Let $p$ be unramified in $K_n$, and let $C_n = \{1,\gamma \} \subset
\gal(K_n/\Q)$, where $\gamma$ is given by 
$\gamma(\zeta_n) = \zeta_n^{-1}$ and $\gamma(\alpha_n) = \alpha_n^{-1}$.
Then the condition that $n|i_p$ is equivalent to 
$\sigma_p \in C_n$.  Moreover, $C_n$ is invariant under conjugation.
\end{lem}

\begin{proof}
{\em The split case:} Since $n|i_p$ and $i_p|p-1$ we have $\zeta_n \in
F_p$, i.e.  $F_p$ contains all $n$-th roots of unity.  Moreover, 
$\unitp$ is an $n$-th power of some element in $F_p$, and thus the
equation $x^n -\unit$ splits completely in $F_p$. In other words, 
$p$ splits completely in $K_n$ and $\sigma_p$ is trivial. 

{\em The inert case:} 
Since $n$ divides $i_p$, $\unitp$ is an $n$-th power of some element
in $F_{p^2}^1$ and hence $\alpha_n \in F_{p^2}$. Moreover,
$n|p^2-1$ implies that $\zeta_n \in F_{p^2}$. Now,
$N_{F_p}^{F_{p^2}}(\alpha_n)=1$ and
$N_{F_p}^{F_{p^2}}(\zeta_n)=\zeta_n^{p+1}=1$ implies that 
$$
\sigma_p(\zeta_n) \equiv \zeta_n^{-1} \mod p, \quad
\sigma_p(\alpha_n) \equiv \alpha_n^{-1} \mod p.
$$
For $p$ that does not ramify in $K_n$ we thus have
\begin{equation}
\label{e:frob-inert}
\sigma_p(\zeta_n) = \zeta_n^{-1}, \quad
\sigma_p(\alpha_n) = \alpha_n^{-1}
\end{equation}
Now, an
element $\tau \in \gal(K_n/\Q)$ is of the form
$$
\tau \colon 
\begin{cases}
\zeta_n \rightarrow \zeta_n^t & t \in \Z \\
\alpha_n \rightarrow \alpha_n^{u} \zeta_n^{s} & 
s \in \Z, \quad u \in \{ 1, -1 \}
\end{cases}
$$
Composing $\gamma$ and $\tau$ then gives
$$
\tau \circ \gamma: 
\begin{cases}
\zeta_n \rightarrow \zeta_n^{-1} \rightarrow \zeta_n^{-t} \\
\alpha_n \rightarrow \alpha_n^{-1} \rightarrow 
\alpha_n^{-u} \zeta_n^{-s} 
\end{cases}
$$
and
$$
\gamma \circ \tau : 
\begin{cases}
\zeta_n \rightarrow \zeta_n^t \rightarrow \zeta_n^{-t} \\
\alpha_n \rightarrow \alpha_n^{u} \zeta_n^{s} \rightarrow 
\alpha_n^{-u} \zeta_n^{-s}
\end{cases}
$$
which shows that $\gamma$ is invariant under conjugation. 
\end{proof}

\subsection{The Chebotarev density Theorem}
In \cite{serre-chebotarev} Serre proved that the Generalized Riemann
Hypothesis (GRH) implies the following version of the
Chebotarev density Theorem: 
\begin{thm}
\label{t:serre-effective-grh}
  Let $E/\Q$ be a finite Galois extension of degree $[E:\Q]$ and
  discriminant $D_E$. For $p$ a prime let $\sigma_p \in G = \gal(E/\Q)$
  denote the Frobenius conjugacy class, and let $C \subset G$
  be a union of conjugacy classes. If the nontrivial zeroes of
  $\zeta_E(s)$ lie on the line $Re(s)=1/2$, then for $x \geq 2$,
$$
|\{ p \leq x \colon \sigma_p \in C \}| 
=
\frac{|C|}{|G|}\pi(x) + 
O\left( 
  \frac{|C|}{|G|}x^{1/2} \left(\log D_E + [E:\Q] \log x \right)
\right)
$$
\end{thm}
Now, primes that ramify in $K_n$ divides $n D_K$ (see
Lemma~\ref{l:discriminant-bound}), so as far as densities are
concerned, ramified primes can be ignored. The bounds on the size of
$D_{K_n}$ (see Lemma~\ref{l:discriminant-bound}) and
Lemma~\ref{l:size-of-frob-conj-class} then gives the following:
\begin{cor}
\label{c:grh-density}
If GRH is true then
\begin{equation}
|\{ p \leq x \colon n|i_p \}| =
 \frac{2}{[K_n:\Q]} \times \pi(x)  + 
O\left( x^{1/2}(\log(xn)  \right)
\end{equation}
\end{cor}

{\em Remark:}  For theorems ~\ref{t:hooley++} and \ref{t:all-N} to be
true, it is enough to assume that the 
Riemann hypothesis holds for all $\zeta_{K_n}$,  $n >1$.

\subsubsection{Bounds on degrees}

In order to apply the Chebotarev density Theorem we need bounds
on the degree $[K_n:\Q]$. We will first assume that
$\unit$ is a fundamental unit. 

\begin{lem} 
\label{l:nonabelian} If $\unit$ is a fundamental unit in $K$ and if
$n=4$ or $n=q$, for $q$ an odd prime, then $\gal(K_n/K)$ is 
  nonabelian.
\end{lem}
\begin{proof} 
  We start by showing that $[K_n:Z_n]=n$. 
  Consider first the case $n=q$. If $\alpha_q \in Z_q$ then $\beta =
  N^{Z_q}_{K}(\alpha_q) = \alpha_q^{[Z_q:K]} \zeta_q^t \in K \subset
  \R$ for some integer $t$.  Since $q$ is odd we may assume that
  $\alpha_q \in \R$, 
  and this forces $\zeta_q^t = 1$, which in turn implies that
  $\alpha_q^{[Z_q:K]} \in K$. Because $\unit$ 
  is a fundamental unit this means that $q|[Z_q:K]$. On the other
  hand, $[Z_q:K]|\phi(q)$, a contradiction. Thus $\alpha_q \not \in
  Z_q$, and hence $K_q/Z_q$ is a Kummer extension of degree $q$.
  
  For $n=4$ we note that $i \in Z_4=K(i)$. Thus $\alpha_2 =
  \sqrt{\unit} \in Z_4$ implies that $\sqrt{-\unit} \in Z_4$. However,
  either $\sqrt{\unit}$ or $\sqrt{-\unit}$ is real and generates a
  {\em real} degree two extension of $K$, whereas $K(i)$ is a
  non-real quadratic extension of $K$, and hence $\alpha_2 \not \in Z_4$.
  Now, if $\alpha_4 \in Z_4(\alpha_2)$ then
  $N^{Z_4(\alpha_2)}_{Z_4}(\alpha_4)= \alpha_4^2 i^t \in Z_4$ for some
  $t\in \Z$, and thus $\alpha_4^2=\alpha_2 \in Z_4$ which contradicts
  $\alpha_2 \not \in Z_4$. Therefore,
$$
[Z_4(\alpha_4):Z_4]
= 
[Z_4(\alpha_4):Z_4(\alpha_2)] [Z_4(\alpha_2):Z_4] =4.
$$ 

Finally we note that the commutator of any nontrivial element
$\sigma_1 \in \gal(K_n/Z_n)$ with any nontrivial element $\sigma_2 \in
\gal(K_n/L_n)$ is nontrivial (we may regard $\gal(K_n/Z_n)$ and
$\gal(K_n/L_n)$ as subgroups of $\gal(K_n/K)$). Hence $\gal(K_n/K)$ is
nonabelian.
\end{proof}

\begin{lem} 
\label{l:degree-lower-bound}
  If $\unit$ is a fundamental unit then
  $$
  [K_n:Z_n] \geq n/2.
  $$
\end{lem}
\begin{proof} 
Clearly $Z_n(\alpha_{q^k}) \subset K_n$, and since field extensions of
relative prime degrees are disjoint,  it is enough to show that if
  $q^k||n$ is a prime power then $q^k|[Z_n(\alpha_{q^k}):Z_n]$ if $q$
  is odd, and $q^{k-1}|[Z_n(\alpha_{q^k}):Z_n]$ if $q=2$.
  
  If $q$ is odd then Lemma~\ref{l:nonabelian} implies that $\alpha_q
  \not \in Z_n$ since $\gal(Z_n/K)$ is abelian. Hence, if $m \in \Z$
  and $\alpha_{q^k}^m \in Z_n$, we must have $q^k|m$. Now, if $\sigma
  \in \gal(Z_n(\alpha_{q^k})/{Z_n})$ then
  $\sigma(\alpha_{q^k})=\alpha_{q^k} \zeta_{q^k}^{t_\sigma}$ for some
  integer $t_\sigma$. Thus there exists an integer $t$ such that
  $$
  \beta = N^{Z_n(\alpha_{q^k})}_{Z_n}(\alpha_{q^k}) =
  \alpha_{q^k}^{[Z_n(\alpha_{q^k}):Z_n]} \zeta_q^t \in Z_n 
  $$
  Multiplying $\beta$ by $\zeta_q^{-t} \in Z_n$ we find that
  $\alpha_{q^k}^{[Z_n(\alpha_{q^k}):Z_n]} \in Z_n$, and hence
  $q^k|[Z_n(\alpha_{q^k}):Z_n]$. 

  For $q=2$ the proof is similar, except that a factor of two is lost
  if $\alpha_2 \in Z_n$.
\end{proof}

{\em Remark:} $K_2/Q$ is a Galois extension of degree four, hence
abelian and therefore contained in some cyclotomic extension by the
Kronecker-Weber Theorem, and it is thus possible that $\alpha_2 \in Z_n$
for some values of $n$.

\begin{lem} 
\label{l:kummer-degree}
We have 
$$
n \phi(n) \ll_K [K_n:\Q] \leq 2 n \phi(n)
$$
\end{lem}
\begin{proof} 
  We first note that $[Z_n:K]$ equals $\phi(n)$ or $\phi(n)/2$
  depending on whether $K \subset \Q(\zeta_n)$ or not. We also
  have the trivial upper bound $[K_n:Z_n] \leq n$. 
  
  For a lower bound of $[K_n:Z_n]$ we argue as follows: Let $\gamma
  \in K$ be a fundamental unit. Since the norm of $\unit$ is one we
  may write $\unit=\gamma^k$ for some $k \in \Z$.  (Note that $k$ does
  not depend on $n$.)  As $[Z_n(\gamma^{1/n}):Z_n(\unit^{1/n})] \leq
  k$, Lemma~\ref{l:degree-lower-bound} gives that
  $[Z_n(\unit^{1/n}):Z_n] \geq n/k$.  The upper and lower bounds now
  follows from
$$
[K_n:\Q] =
[K_n:Z_n][Z_n:K][K:\Q]
$$
\end{proof}

\subsubsection{Bounds on discriminants}

\begin{lem}
\label{l:discriminant-bound}
If $p$ ramifies in $K_n$ then $p|nD_K$. Moreover,  
$$
\log( \disc(K_n/\Q) ) \ll_K [K_n:K] \log(n)
$$
\end{lem}
\begin{proof}
First note that 
$$
\disc(K_n/\Q) = N^K_{\Q}(\disc(K_n/K)) \times \disc(K/\Q)^{[K_n:K]}.
$$
From the multiplicativity of the different we get
$$
\disc(K_n/K) = \disc(Z_n/K)^{[K_n:Z_n]} \times N^{Z_n}_K( \disc(K_n/Z_n)),
$$
Since $\unit$ is a unit, so is $\unit^{1/n}$. Thus, if we let
$f(x)=x^n-\unit$ then $f'(x)=nx^{n-1}$, and therefore the principal
ideal $f'(\unit^{1/n})
\OO_{K_n}$ equals $n \OO_{K_n}$. In terms of discriminants this means that 
$$
\disc(K_n/Z_n) | N^{K_n}_{Z_n}( n \OO_{K_n} )
$$
and similarly it can be shown that 
$$
\disc(Z_n/K) | N^{Z_n}_K( n \OO_{Z_n} ).
$$
Thus $\disc(K_n/\Q)$ divides
$$
N^K_{\Q} \left( 
  N^{K_n}_{K}( n \OO_{K_n} ) \times 
  N^{Z_n}_K( n \OO_{Z_n} )^{[K_n:Z_n]}
\right) 
\times 
\disc(K/\Q)^{[K_n:K]}.
$$
$$
=
n^{ 4[K_n:K]} \times
\disc(K/\Q)^{[K_n:K]}
$$
which proves the two assertions. 
\end{proof}
\section{Proof of Theorem \ref{t:hooley++}}

In order to bound the number of primes $p<x$ for which $i_p>x^{1/2}$
we will need the following Lemma:
\begin{lem} 
\label{l:simple}
The number of primes $p$ such that $\ord_p(A) \leq y$ is $O(y^2)$. 
\end{lem}
\begin{proof} 
  Given $A$ there exists a constant $C_A$ such that $\det(A^n-I) =
  O(C_A^n)$. Now, if the order of $A$ mod $p$ is $n$, then certainly
  $p$ divides $\det(A^n-I) \neq 0$. Putting $M=\prod_{n=1}^{y}
  \det(A^n-I)$ we see that any prime $p$ for which $A$ has order $n
  \leq y$ must divide $M$. Finally, the number of prime divisors of
  $M$ is bounded by
$$
\log(M) \ll \sum_{n=1}^{y} n \log(C_A) \ll y^2.
$$
\end{proof}

{\em First step:} We consider primes $p$ such that $i_p \in (x^{1/2}
\log x, x)$. By Lemma~\ref{l:simple} the number of such primes is
\begin{equation}
O \left( \left( \frac{ x }{ x^{1/2} \log x } \right)^2 \right) =
O \left( \frac{ x }{ \log^2 x } \right).
\end{equation}

{\em Second step:} Consider $p$ such that $q | i_p$ for some prime $q
\in ( \frac{x^{1/2}}{\log^3 x}, x^{1/2} \log x )$. We may bound this
by considering primes $p \leq x$ such that $p \equiv \pm 1 \mod q$ for
$q \in ( \frac{x^{1/2}}{\log^3 x}, x^{1/2} \log x )$. Since $q \leq
x^{1/2}\log x $, Brun's sieve gives (up to an absolute constant)
the bound
$$
\frac{x}{\phi(q) \log(x)}
$$
and the total contribution from these primes is at most
\begin{equation}
\label{e:small-range}
\sum_{ q \in ( \frac{x^{1/2}}{\log^3 x}, x^{1/2} \log x )  }
\frac{x}{\phi(q) \log(x/q)}
\ll
\frac{x}{\log x}
\sum_{ q \in ( \frac{x^{1/2}}{\log^3 x}, x^{1/2} \log x )  }
\frac{1}{q}.
\end{equation}
Now, summing reciprocals of primes in a dyadic interval, we get
$$\sum_{q \in [M,2M]} \frac{1}{q} \ll \frac{ \pi(2M)}{M} \leq
\frac{1}{\log M}$$
Hence 
$$
\sum_{ q \in ( \frac{x^{1/2}}{\log^3 x}, x^{1/2} \log x )  }
\frac{1}{q} 
\ll 
\frac{1}{\log x} 
\log_2 \left( \frac{ x^{1/2} \log x  }{x^{1/2}/\log^3 x} \right)
\ll 
\frac{ \log \log x}{\log x}.
$$
and equation~\eqref{e:small-range} is $O( \frac{x \log\log
  x}{\log^2 x})$. 

{\em Third step:} Now consider $p$ such that $q| i_p$ for some prime $q \in
( f(x)^2, \frac{x^{1/2}}{\log^3 x})$. We are now in the range where
GRH is applicable; by Corollary~\ref{c:grh-density} and
Lemma~\ref{l:kummer-degree} we have 
$$
|\{ p \leq x : q | i_p  \}|
\ll 
\frac{x}{q \phi(q) \log x} + O( x^{1/2} \log( x q^2) )
$$
Summing over $q \in (f(x)^2,\frac{x^{1/2}}{\log^3 x})$ we find that the 
number of such $p \leq x$  is bounded by 
\begin{equation}
\label{e:grh-bound-two}
\sum_{ q \in ( f(x)^2, \frac{x^{1/2}}{\log^3 x})  }
\left( 
  \frac{x}{q^{2} \log x} + O( x^{1/2} \log( x q^2) )
\right)
\end{equation}
Now,
$$
\sum_{ q \in ( f(x)^2, \frac{x^{1/2}}{\log^3 x})  }
\frac{1}{q^{2}}
\ll 
\frac{1}{f(x)}
$$
and thus equation \eqref{e:grh-bound-two} is 
$$
\ll
\frac{x}{f(x) \log x}
+
\frac{x}{\log^2 x}.
$$

{\em Fourth step:} For the remaining primes $p$, any prime divisor
$q|i_p$ is smaller than $f(x)^2$. Hence $i_p$ must be divisible by
some integer $d \in (f(x),f(x)^3)$. Again Lemmas~\ref{c:grh-density}
and \ref{l:kummer-degree} give
$$
|\{ p \leq x : d | i_p  \}|
\ll 
\frac{x}{d \phi(d) \log x} + O( x^{1/2} \log( x d^2) )
$$
Noting that $\phi(d)\gg d^{1-\epsilon}$ and summing over $d \in (f(x),
f(x)^3)$ we find that the 
number of such $p \leq x$  is bounded by 
\begin{equation}
\label{e:grh-bound-one}
\sum_{ d \in (f(x), f(x)^3)  }
\left( 
  \frac{x}{d^{2-\epsilon} \log x} + O( x^{1/2} \log( x d^2) )
\right)
\end{equation}
Now,
$$
\sum_{ d \in ( f(x),f(x)^3)  }
\frac{1}{d^{2-\epsilon}}
\ll 
\frac{1}{f(x)^{1-\epsilon}}
$$
and
$$
\sum_{ d \in ( f(x),f(x)^3)  }
x^{1/2} \log( x d^2) 
\ll
f(x)^3 x^{1/2} \log( x^2) 
$$
therefore equation \eqref{e:grh-bound-one} is 
$$
\ll
\frac{x}{f(x)^{1-\epsilon} \log x}
$$


\section{Proof of Theorems \ref{t:ord-b-all-N} and \ref{t:all-N}}

Given a composite integer $N = \prod_{p|N} p^{a_p}$ we
wish to use the lower bounds on $\ord_p(b)$ (or $\ord_p(A)$) to obtain
a lower bound on $\ord_N(b)$.  The main obstacle is that $\ord_N(b)$
can be much smaller than $\prod_{p|N} \ord_{p^{a_p}}(b)$.  Let
$\lambda(N)$ be the Carmichael lambda function, i.e., the exponent of
the multiplicative group $(\Z/N\Z)^\times$.  Clearly $\ord_N(b) \leq
\lambda(N)$, and it turns out that $\lambda(N)$ can be much smaller
than $N$.
However, 
 $\lambda(N)\gg N^{1-\epsilon}$ for most
$N$
(see \cite{carmichael-lambda}), and since
$$
\ord_N(b) \geq 
\frac{\lambda(N)}{N}
\prod_{p|N} \ord_p(b)
$$
it suffices to show that most integers are essentially given by a
product of primes $p$ such that $\ord_p(b) \geq p/\log p$.  
%
%
We will
only give the details for Theorem \ref{t:all-N} since the other case
is very similar.

If $p$ is prime such that $\ord_p(A) \leq p/\log(p)$, or $p$ ramifies
in $K$, we say that $p$ is
``bad''. We let $P_B$ denote the set of all bad primes, and we let
$P_B(z)$ be the set of primes $p \in P_B$ such that $p \geq z$. 
Since only finitely many primes ramify in $K$, Theorem
\ref{t:hooley++} gives that  the number of bad primes $p \leq x$ is
$O(\frac{x}{\log^{2-\epsilon} x})$. A key observation is the following: 
\begin{lem} 
We have
\begin{equation}
\label{e:convergence}
\sum_{p \in P_B} \frac{1}{p} < \infty
\end{equation}
In particular, if we let $$\beta(z) = \sum_{p \in P_B(z)} 1/p,$$ then
$\beta(z)$ tends to zero as  $z$ tends to infinity.
\end{lem}
\begin{proof}
Immediate from partial summation and the $O(\frac{x
  }{\log^{2-\epsilon} x})$ estimate in Theorem~\ref{t:hooley++}.
\end{proof}

Given $N \in Z$, write $N = s^2 N_G N_B$ where $N_G N_B$ is square
free and $N_B$ is 
the product of ``bad'' primes dividing $N$.  By the following Lemma, we
find that  few integers have a large square factor:
\begin{lem} We have
$$
| \{ N \leq x \ : \ s^2|N, s \geq y \}| = 
O \left(\frac{x}{y} \right)
$$
\end{lem}
\begin{proof}
The number of $N \leq x$ such that $s^2 |N$ for $s \geq y $ is
bounded by
$\sum_{s \geq y} \frac{x}{s^2} \ll \frac{x}{y}.$
\end{proof}

Next we show that there are few $N$ for which $N_B$ is divisible by $p
\in P_B(z)$.  In other words, for most $N$, $N_B$ is a product of
small ``bad'' primes.

\begin{lem}
The number of $N \leq x$ such that $p \in P_B(z)$ divides $N_B$ is
$O(x \beta(z))$. 
\end{lem}
\begin{proof}
Let $p \in P_B(z)$. The number of $N \leq x$ such that $p|N$ is less
than $x/p$. Thus, the total number of $N \leq x$ such that some $p\in
P_B(z)$ divides $N$, is bounded by
$$
\sum_{p \in P_B(z)} \frac{x}{p} = x \sum_{p \in P_B(z)} \frac{1}{p} = x \beta(z).
$$
\end{proof}

Combining the previous results we get that the number of $N = s^2 N_G
N_B \leq x$ 
such that $N_B$ is $z$-smooth and $s \leq y$ is
$$
x \left( 1 + O\left( \beta(z)+1/y \right) \right).
$$
For such $N$ we have $N_B \leq \prod_{p \leq z} p \ll e^z$. Letting
$z=\log \log x$ and $y=\log x$ we get that 
$$
N_G = \frac{N}{s^2 N_B} \geq \frac{N}{\log^3 x}
$$
for $N \leq x$ with at most $O\left( x( \beta(\log \log x) + (\log
  x)^{-1}) \right) = o(x)$ 
exceptions. 
Now, the following Proposition gives that, for most $N$, $\ord_N(A)$ is
essentially given by $\prod_{p|N} \ord_p(A)$.
\begin{prop2}[\cite{cat2}, Proposition 11]
\label{prop L}
Let $D_A = 4(\trace(A)^2-4)$.  For almost all
\footnote{By ``for almost all $N \leq x$'' 
we mean that there are $o(x)$ exceptional integers $N$ that are
smaller than $x$.} 
$N\leq x$,
$$
\ord_N(A) \geq \frac{\prod_{p\mid d_0} \ord_p(A)}
{ \exp(3(\log\log x)^4)}
$$
where $d_0$ is given by writing $N=ds^2$, with $d=d_0\gcd(d,D_A)$
square-free.
\end{prop2}
Finally, since $\ord_p(\A) \geq \frac{p}{\log p} \geq
p^{1-\epsilon}$ for $p|N_G$ and $p$ sufficiently large, we find
that 
$$
\ord_N(\A) \gg 
\frac{\prod_{p\mid N_G} \ord_p(A)}
{ \exp(3(\log\log x)^4)}
\gg
\frac{N_G^{1-\epsilon}} 
{ \exp(3(\log\log x)^4)}
\gg N^{1-2\epsilon}
$$
for all but $o(x)$ integers $N \leq x$.

\bibliographystyle{abbrv} 

\end{document}